\documentclass[a4paper,12pt]{article}
\usepackage{amsmath,amssymb}
\usepackage{amscd}

\usepackage[usenames]{color}
\usepackage[pdfpagelabels, pdftex]{hyperref}

\title{Descent obstruction and fundamental exact sequence}

\author{David Harari and Jakob Stix}
\date{May 26, 2010}

\DeclareFontFamily{U}{wncy}{}
\DeclareFontShape{U}{wncy}{m}{n}{%
   <5>wncyr5%
   <6>wncyr6%
   <7>wncyr7%
   <8>wncyr8%
   <9>wncyr9%
   <10>wncyr10%
   <11>wncyr10%
   <12>wncyr6%
   <14>wncyr7%
   <17>wncyr8%
   <20>wncyr10%
   <25>wncyr10}{}
\DeclareMathAlphabet{\cyrille}{U}{wncy}{m}{n}
\def\Sha{\cyrille X}

\def\Ga{\Gamma}

\def \A{\bf A}

\def \rH{{\rm H}}
\def \out{{\rm out}}

\newcommand{\surj}{\twoheadrightarrow} 

\DeclareMathOperator{\Spec}{Spec}

\newcommand{\dO}{{\mathcal O}}

\newcommand{\bC}{{\mathbb C}}
\newcommand{\bP}{{\mathbb P}}
\newcommand{\bQ}{{\mathbb Q}}
\newcommand{\bR}{{\mathbb R}}
\newcommand{\bZ}{{\mathbb Z}}

\newtheorem{theo}{Theorem}[section]
\newtheorem{prop}[theo]{Proposition}
\newtheorem{lem}[theo]{Lemma}
\newtheorem{cor}[theo]{Corollary}

\def \kbar {{\bar k}}

\def \dem {\paragraph{Proof:}}

\def \rem {\paragraph{Remark:}}
\def \rems {\paragraph{Remarks:}}

\def \Romannumeral #1 {\expandafter\uppercase\expandafter {\romannumeral #1} }
\def \br {{\rm{Br\,}}}

\def \gal {{\rm{Gal\,}}}

\def\ov{\overline}

\def \Hom {{\rm {Hom}}}

\def \Q {{\bf Q}}

\def\smallsquare{\vbox{\hrule\hbox{\vrule height 1 ex\kern 1 ex\vrule}\hrule}}
\def\enddem{\hfill \smallsquare\vskip 3mm}
\def \merci {\paragraph{Aknowledgements.}}
\def \abstract{\paragraph{Abstract.}}

\newcommand{\df}{\bf}   

\hypersetup{
  pdftitle={},
  pdfauthor={},
  pdfsubject={},
  pdfkeywords={},
  colorlinks=true,
  breaklinks=true,
  bookmarksopen=true,
  bookmarksnumbered=true,
  pdfpagemode=UseOutlines,
  plainpages=false}

\begin{document}
\maketitle


\subsection*{Introduction}

Let $k$ be a field of characteristic zero 
with algebraic closure $\kbar$ and absolute Galois group 
$\Ga_k:=\gal(\kbar/k)$. Let $X$ be 
a geometrically connected variety over $k$. Fix a 
geometric point $\bar x \in X(\kbar)$ and let $\pi_1(X):=\pi_1(X,\bar x)$
be the \'etale fundamental group of $X$. Set $\ov X=
X \times_k \kbar$ and denote by $\pi_1(\ov X):=\pi_1(\ov X,\bar x)$
the \'etale fundamental group of $\ov X$.
Recall Grothendieck's fundamental exact sequence of profinite groups, cf.\ \cite{sga1} IX Thm 6.1,
\begin{equation} \label{fund}
1 \to \pi_1(\ov X) \to \pi_1(X) \to \Ga_k \to 1.
\end{equation}
By covariant functoriality of $\pi_1$, the existence of a $k$-point on $X$ 
implies that the exact sequence (\ref{fund}) has a section. 
Grothendieck's section conjecture predicts that the converse statement 
is true whenever $X$ is a proper\footnote{There is also a version of the section conjecture for affine hyperbolic curves when $k$-rational cusps need to be considered as well, see \cite{letter} page 8/9.} hyperbolic curve over a number field, see \cite{letter}, (and there is 
also a $p$-adic version of this conjecture).

The goal of this note is to relate the existence of a section for 
(\ref{fund}) when $k$ is a number field to the fact that $X$ has an adelic point for which 
there is no descent obstruction (in the sense of \cite{skobook}  II section 5.3) 
associated to torsors under finite group 
schemes, see Theorem~\ref{numbers}. 

We will also prove related statements over arbitrary fields, so that the reader 
can distinguish between purely formal results and results related 
to arithmetic properties. 

\medskip

Since we want to deal with variants of the exact sequence (\ref{fund}), for example the
abelianized fundamental exact sequence, we need 
to introduce a general setting as follows.
Let $\ov U < \pi_1(\ov X)$ be a closed subgroup  which is a normal subgroup in  $\pi_1(X)$.
For example $\ov U$ could be a characteristic subgroup in $\pi_1(\ov X)$ of which there are plenty, because $\pi_1(\ov X)$ is finitely generated as a profinite group.
We set $\ov A=\pi_1(\ov X)/\ov U$. 
The pushout of (\ref{fund}) by the canonical surjection $\pi_1(\ov X) \to \ov A$ is the exact 
sequence
\begin{equation} \label{modfund}
1 \to \ov A \to \pi_1 ^U (X) \to \Ga_k \to 1,
\end{equation}
that can be defined because $\ov U$ is also normal in $\pi_1(X)$ and in particular by definition the kernel of the induced quotient map $\pi_1(X) \to \pi_1 ^U(X)$.
This construction contains as special cases for $\ov A$ the 
profinite abelianized group $\pi_1(\ov X)^{\rm ab}$ if 
we take for $\ov U$ the closure of the derived 
subgroup of $\pi_1(\ov X)$, and $\ov A = \pi_1(\ov X)$ if we  take for $\ov U$ the trivial group.

\smallskip

Let $(k_i)_{i \in I}$ be a family of field extensions of $k$. The guiding example is the case of a number field $k$ together with the family of all completions $k_v$ of $k$, 
or of all corresponding henselizations\footnote{By convention if 
$k_v$ is an archimedean completion of $k$, the henselization 
$k_v^h$ means the algebraic closure of $k$ into $k_v$.} $k_v^h$ of $k$. Fix
an algebraic closure $\kbar_i$ of $k_i$ and embeddings $\kbar \to 
\kbar_i$, so that with $\Ga_i :=\gal(\kbar_i/k_i)$ we have canonical restriction maps $\theta_i : \Ga_i \to \Ga_k$.
Moreover,  the embedding $\kbar \to \kbar_i$ yields canonically a geometric point $\bar x_i$ of $X_i:= X \times_k k_i$, that projects onto $\bar x$, and thus a canonical map $\pi_1(X_i, \bar x_i) \to \pi_1(X, \bar x)$. The reason for assuming characteristic zero in the first place is that by the 
comparison theorem (cf. \cite{szapi}, p. 186, Remark~5.7.8.),
the natural maps induce an isomorphism
\[
\pi_1(X_i,\bar x_i) \xrightarrow{\sim} \pi_1(X,\bar x) \times_{\Ga_k} \Ga_i.
\]
Thus a section $s$ of (\ref{fund}) induces canonically a section $s_i$ of the analogue of (\ref{fund}) for the $k_i$-variety $X_i$. In this formal setting the main goal of this note is to establish a criterion inspired by the descent obstruction for when a collection
of sections $(s_i)_{i \in I}$ comes from a section $s$ up to conjugation from $\pi_1(\ov X)$.


\subsection*{Reminder on nonabelian $\rH^1$}

Recall that if $G$ is a 
finite $k$-group scheme, then 
the \'etale cohomology set $\rH^1(X,G)$ is the same as the
Galois cohomology set $\rH^1(\pi_1(X),G(\kbar))$, where the action 
of $\pi_1(X)$ on $G(\kbar)$ is induced by the projection map $\pi_1(X) \to \Ga_k$ that occurs in (\ref{fund}), cf.\ \cite{sga1} XI \S 5.

The identification is natural in both $X$ and $G$, although for a map $Y \to X$ the induced map $\pi_1(Y) \to \pi_1(X)$ is only well defined up to inner automorphism by an element of $\pi_1(\ov X)$. In fact, such an inner automorphism acts as the identity on $\rH^1(\pi_1(X),G)$ by the following reasoning.
Recall that in general if 
$\varphi = \gamma(-)\gamma^{-1}$ is an inner automorphism of a profinite 
group $\pi$ by an element $\gamma \in \pi$ and
$M$ is a discrete $\pi$-group, then the 
map $(\varphi^*,\gamma^{-1}) : \rH^1(\pi,M) \to \rH^1(\pi,M)$ which is the composite
\[
\rH^1(\pi,M) \xrightarrow{\varphi^\ast} \rH^1(\pi,\varphi^\ast M) \xrightarrow{\gamma^{-1}.} \rH^1(\pi,M),
\]
that exploits the $\pi$-map "multiplication by $\gamma^{-1}$" $: \varphi^\ast M \to M$ and 
which on cocycles is given by 
\[
\big(\sigma \mapsto a_\sigma\big) \mapsto \big(\sigma \mapsto \gamma^{-1}(a_{\gamma \sigma\gamma^{-1}})\big)
\]
is the identity map. This is classical when $M$ is abelian, see \cite{col}, VII.5.\ Proposition 3, 
and easy to check in the general case by the same direct computation:
\begin{equation} \label{conjugation}
\gamma^{-1}(a_{\gamma \sigma\gamma^{-1}}) = \gamma^{-1}(a_\gamma) a_{\sigma \gamma^{-1}} = 
 \gamma^{-1}(a_\gamma) a_\sigma \sigma(\gamma^{-1}(a_\gamma)),
\end{equation}
which shows that $(\sigma \mapsto a_{\sigma})$ is indeed cohomologous 
to $\sigma \mapsto \gamma^{-1}(a_{\gamma \sigma\gamma^{-1}})$.
In our geometric example the element $\gamma \in \pi_1(\ov X)$ acts trivially on the coefficients $M = G(\kbar)$ such that  $(\varphi^\ast,\gamma^{-1})$ becomes simply the pullback by conjugation with $\gamma$, which therefore acts as identity on $\rH^1(\pi_1(X),G(\kbar))$.

The \'etale cohomology set 
$\rH^1(\ov X,\ov G)$  for $\ov G:=G \times_k \kbar$ is naturally the set 
$\Hom^\out(\pi_1(\ov X),G(\kbar))$ 
of continuous
homomorphisms $\pi_1(\ov X) \to G(\kbar)$ up to conjugation
by an element of $G(\kbar)$, see \cite{cogal} I.5.

\medskip

The following well known interpretation of $\rH^1(\Ga,G)$ will become useful later. Let 
\[
1 \to G \to E \to \bar \Ga \to 1
\]
be short exact sequence of profinite groups, and let 
$\ov \varphi : \Ga \to \bar \Ga$ 
be a continuous homomorphism. The set of lifts $\varphi : \Ga \to E$ of $\ov \varphi$ up to conjugation by an element of $G$ is either empty or, with the group $G$ equipped with the conjugation action of $\Ga$ via a choice of lift $\varphi_0$, in bijection with the corresponding $\rH^1(\Ga,G)$. Indeed for a cocycle $a : \Ga \to G$ the twist of $\varphi_0$ by $a = (\gamma \mapsto a_\gamma)$, i.e. the map 
\[
\gamma \mapsto a_\gamma \cdot \varphi_0(\gamma),
\] 
is another lift. Any other lift of $\ov \varphi$ can be described by such a twist. Two cocycles are cohomologous if and only if they lead to conjugate lifts. The description of lifts via $\rH^1(\Ga,G)$ is natural with respect to both $\Ga$ and $G$.


\section{Results over arbitrary fields} 

The notation and assumptions in this whole section are as above.
In particular we consider the exact sequence (\ref{modfund})
associated to a quotient $\ov A$ of $\pi_1(\ov X)$ by a  
subgroup $\ov U$ that remains normal in $\pi_1(X)$. 

\smallskip

For each $i \in I$, let $\sigma_i : \Ga_i \to \pi_1(X_i,\bar{x}_i)$ be a section of the fundamental sequence associated to $X_i$. For example $\sigma_i$ could be the section associated to a $k_i$-point $P_i \in X(k_i)$.
By composition we obtain a {\df section map}
\[
s_i : \Ga_i \xrightarrow{\sigma_i} \pi_1(X_i) \to \pi_1(X) \to \pi_1 ^U(X). 
\]

Let $G$ be a finite $k$-group scheme, hence $G(\kbar)$ is a finite discrete $\Ga_k$-group. Via $\theta_i:\Ga_i \to \Ga_k$ we may view $G(\kbar) = \theta_i^\ast G(\kbar)$ also as a discrete $\Ga_i$-group that describes the base change $G \times_k k_i$. 

A cohomology class $\alpha \in \rH^1(X,G)$  such that the corresponding element $\bar \alpha \in 
\rH^1(\ov X,\ov G)$ has trivial restriction to $\ov U$, has an evaluation $\alpha(s_i) \in \rH^1(k_i,G)$ as follows. By the restriction--inflation sequence the class $\alpha$ uniquely comes from $\rH^1(\pi_1^U(X),G(\kbar))$ and so the pullback class $\alpha(s_i) := s_i^\ast(\alpha) \in \rH^1(k_i,G)$ is defined. 
Note that the coefficients $G$ here are indeed the group $G(\kbar)$ with $\Ga_i$ action induced by $\theta_i$ because $s_i$ comes from a section $\sigma_i$. By formula (\ref{conjugation}) the evaluation does only depend on $s_i$ up to conjugation by an element of $\pi_1^U(X)$ with trivial action on $G(\kbar)$.

By analogy with \cite{stollant}, Definition~5.2.,
we say that the tuple of section maps $(s_i)_{i \in I}$ {\df survives every finite descent obstruction} if the following holds.

\begin{enumerate}
\item[(a)] 
For every finite $k$-group scheme $G$ and every $\alpha \in 
\rH^1(X,G)$ such that the corresponding element $\bar \alpha \in 
\rH^1(\ov X,\ov G)$ has trivial restriction to $\ov U$, the family 
$(\alpha(s_i)) \in \prod_{i \in I} 
\rH^1(k_i,G)$ belongs to the diagonal image of $\rH^1(k,G)$. 
\end{enumerate}
Clearly, if the sections $s_i$ are the sections associated to $k_i$-rational points, then $(s_i)$ survives every finite descent obstruction if and only if the collection $(P_i)$ of rational points survives every finite descent obstruction in the sense of \cite{stollant}.
We furthermore say that the tuple of section maps $(s_i)_{i \in I}$ {\df survives every finite constant descent obstruction} if the following holds.

\begin{enumerate}
\item[(a')]
For every finite constant $k$-group scheme $G$ and every $\alpha \in 
\rH^1(X,G)$ such that the corresponding element $\bar \alpha \in 
\rH^1(\ov X,\ov G)$ has trivial restriction to $\ov U$, the family 
$(\alpha(s_i)) \in \prod_{i \in I} 
\rH^1(k_i,G)$ belongs to the diagonal image of $\rH^1(k,G)$. 
\end{enumerate}

We first establish a link to 
continuous homomorphisms $\Ga_k \to \pi_1 ^U(X)$. 
Let us define a {\df continuous quotient} of a profinite group as a quotient 
by a normal and closed subgroup. 

\begin{prop} \label{converse}
Consider the following assertion: 
\begin{enumerate}
\item[(b)]
There exists a continuous homomorphism 
$s : \Ga_k \to \pi_1 ^U(X)$ such that 
for each $i \in I$, we have $s_i=s \circ \theta_i$ up to conjugation in 
$\pi_1 ^U (X)$. 
\end{enumerate}
Then (b) implies property (a'). If we assume further that the following hypothesis holds: 
\begin{enumerate}
\item[($\ast$)] For every finite and constant $k$-group scheme $G$, the fibres of the diagonal map 
$\rH^1(k,G) \to \prod_{i \in I} \rH^1(k_i,G)$ are finite.
\end{enumerate}
Then (b) is equivalent to property (a').
\end{prop}

\dem Assume (b). Let $G$ and $\alpha \in \rH^1(X,G)=\Hom^\out(\pi_1(X),G(\kbar))$ 
be as in property (a'). Since the restriction of 
$\alpha$ to $\ov U$ is trivial, 
the class $\alpha$ corresponds to a map $\alpha: \pi_1^U(X) \to G(\kbar)$ up to conjugation in $G(\kbar)$. We get  
\[
\alpha(s_i) = \alpha \circ s_i = \alpha \circ s \circ \theta_i = \theta_i^\ast( \alpha(s))
\]
up to conjugation in $G(\kbar)$, so that  $(\alpha(s_i))$ is the image of 
$\alpha(s) \in \rH^1(k,G)$ under the diagonal map, whence property (a'). 

\smallskip

Suppose now that assertion (a') and the additional hypothesis ($\ast$) hold. We are going to show 
that (b) holds as well. For a finite continuous quotient
$p:\pi_1^U(X) \to G$ we consider the set
\[
S_G:=\{ s' \in \Hom(\Ga_k, G) \ ;  \quad \forall i \in I, \quad
\theta_i ^*(s')=p \circ s_i \in \rH^1(k_i,G)\},
\]
where we view $G$ as a constant $k$-group scheme. The set $S_G$ is 
non empty by assumption (a') and finite thanks to ($\ast$) and to the finiteness of $G$.
Therefore $\varprojlim_G S_G$ where $G$ ranges over all finite continuous quotients of $\pi_1^U(X)$ 
is not empty, see \cite{bourbaki:top1-4} Chapter I \S9.6 Proposition 8.
An element $s \in \varprojlim_G S_G$ is nothing but 
a continuous homomorphism
\[
s : \Ga_k \to \varprojlim_G G = \pi_1^U(X)
\]
such that for every 
$i \in I$, and every finite continuous quotient $p: \pi_1^U(X) \to G$ the equality $p \circ s \circ \theta_i = p \circ s_i$ holds up to conjugation by elements from the finite set
\[
C_{i,G} =\{ c \in G, \quad p \circ s \circ \theta_i=c\big(p\circ s_i\big)c^{-1} \} \subset G
\]
The set $\varprojlim_G C_{i,G}$ is not empty by the same argument, which implies that for each $i \in I$, 
we have $s \circ \theta_i=s_i$ up to conjugation in $\pi_1^U(X)$.
\enddem 

\rem  (1) Without additional assumptions, we cannot force the 
supplementary property that $s$ is a section. Indeed, take $k_i=\kbar$ for 
every $i \in I$. Then all sets $\rH^1(k_i,G)$ are trivial, hence 
the condition (a) and thus condition (a') is automatically satisfied. 
Although condition (b) also holds trivially by the choice of the trivial homomorphism, because there is no interpolation property to be satisfied, nevertheless (\ref{modfund}) does not always admit a section, see for example 
\cite{stix} or \cite{dhsza} for counterexamples over local and global fields. 

(2) For an example with a nontrivial homomorphism $s: \Ga_k \to \pi_1(X)$ but no section we consider the case $k=\bR$, $k_i = \bC$ and a real Godeaux--Serre variety. Computations with SAGE show that the homogenous equations
\begin{eqnarray*}
 z_0^2 + z_1^2 + z_2^2 + z_3^2 + z_4^2 + z_5^2 + z_6^2  & = & 0  \\
 z_0 z_2 + z_1 z_3 + z_4^2 + z_5 z_6  & =  & 0  \\
 i  ( z_0^2  - z_1^2) + 3i (z_2^2 - z_3^2) - 2 z_6^2 &  = & 0   \\
 i (z_0 z_1 + z_2 z_3) + z_4 z_5 + z_5 z_6 + z_6 z_4  & = & 0  
\end{eqnarray*}
define a smooth surface $Y$ of general type in $\bP^6_\bC$ with ample canonical bundle $\omega_Y = \dO(1)|_Y$ as computed by the adjunction formula. By the Lefschetz theorem on hyperplane sections $Y$ is simply connected. The surface $Y$ is preserved by the $\Ga_\bR$-semilinear action of $G=\bZ/4\bZ$ on $\bP^6_{\bC}$ generated by 
\[
[z_0:z_1:z_2:z_3:z_4:z_5:z_6] \mapsto [-\bar{z}_1:\bar{z}_0:-\bar{z}_3:\bar{z}_2:\bar{z}_4:\bar{z}_5:\bar{z}_6],
\]
and avoids the fixed point set of the $G$-action.
Hence the quotient map $Y \to X = Y/G$ is the universal cover of the geometrically connected $\bR$-variety $X$ with $\pi_1(X) = \bZ/4\bZ$. The analogue of (\ref{fund}) for $X$ is given by 
\[
1 \to \pi_1(X \times_\bR \bC) \to \bZ/4\bZ  \to \Ga_{\bR} \to 1
\]
which clearly does not admit sections. Nevertheless, there is a nontrivial homomorphism $\Ga_{\bR} \to \pi_1(X)$.

\medskip

However, if we suppose that (\ref{modfund}) has
a section, then we can prove the following stronger approximation result.


\begin{prop} \label{general}
Consider the following assertion: 
\begin{enumerate}
\item[(c)] There exists a section $s : \Ga_k \to \pi_1 ^U(X)$ of (\ref{modfund}) such that
for each $i \in I$, we have $s_i=s \circ \theta_i$ up to conjugation in
$\ov A$.
\end{enumerate}
Then (c) implies (a) which implies the following (a'').
\begin{enumerate}
\item[(a'')] For every finite $k$-group scheme $G$ and every $\alpha \in
\rH^1(X,G)$ such that the corresponding element $\bar \alpha \in
\rH^1(\ov X,\ov G)$ has trivial restriction to $\ov U$ and is  surjective 
(or equivalently, the $G$-torsor $Y \to X$ corresponding to $\alpha$ 
is assumed to be geometrically connected), the family
$(\alpha(s_i)) \in \prod_{i \in I}
\rH^1(k_i,G)$ belongs to the diagonal image of $\rH^1(k,G)$.
\end{enumerate}
If we moreover assume that ($\ast$) holds and that the exact sequence (\ref{modfund}) 
admits a section $s_0$, then the properties (a), (a'') and (c) are all equivalent.
\end{prop}

\dem  The implications (c) $\Rightarrow$ (a) $\Rightarrow$ (a'') are obvious because by formula (\ref{conjugation}) a section $s$ as in (c) implies for every $\alpha \in \rH^1(X,G)$ as in (a) that
\[
\alpha(s_i) = s_i^\ast(\alpha) = (s \circ \theta_i)^\ast(\alpha) = \theta_i^\ast(s^\ast(\alpha)).
\]
It remains to show that 
(a'') $\Rightarrow$ (c) under the additional assumption of ($\ast$) and 
the existence of 
a section $s_0$ of $\pi_1^U(X) \to \Ga_k$. The method is similar to \cite{stollbis},
Lemma 9.13, which deals with the case when $k$ is a number field 
and $(k_i)$ is the family of its completions for miscellaneous $\ov A$, 
like $\ov A=\pi_1(\ov X)$ or $\ov A=\pi_1^{\rm ab}(\ov X)$. For  the
convenience of the reader we give a grouptheoretic version of the argument.

\smallskip

Let us assume (a''). Let $\ov V < \pi_1(\ov X)$ be an open subgroup of finite index containing $\ov U$ and normal in $\pi_1(X)$. Let $\ov A_V$ be the quotient $\pi_1(\ov X)/\ov V$ and let $p_V :\pi_1^U(X) \surj \pi_1^V(X)$ be the corresponding quotient map. The composition $s_{0,V} = p_V \circ s_0$ splits the exact sequence (\ref{modfund}) for $\ov V$
\begin{equation} \label{modfund-V}
1 \to \ov A_V \to \pi_1^V(X) \to \Ga_k \to 1,
\end{equation}
so that $\pi_1^V(X)$ is isomorphic to a semi-direct product. The map $p_V$ and 
\[
p_{0,V} : \pi^U_1(X) \to \Ga_k \xrightarrow{s_{0,V}} \pi_1^V(X)
\]
lift the natural projection $\pi_1^U(X) \to \Ga_k$. Their difference 
$\gamma \mapsto p_V(\gamma)p_{0,V}(\gamma)^{-1}$
 is a cohomology class $\alpha_V \in \rH^1(\pi_1^U(X),\ov A_V)$, with $\pi_1^U(X)$ acting via $p_{0,V}$ and conjugation, that  corresponds to a class in $\rH^1(X,\ov A_V)$ which becomes trivial when restricted to $\ov U$. 
 The restriction of $\alpha_V$ to $\pi_1(\ov X)$  
equals  the surjective map 
$p_V|_{\pi_1(\ov X)} : \pi_1(\ov X) \surj \ov A_V$, hence is 
geometrically connected. 

We now apply (a'') to the class $\alpha_V$. The class $\alpha_V(s_i) = s_i^\ast(\alpha_V)$ measures the difference between $p_V \circ s_i$ and $p_{0,V} \circ s_i =  s_{0,V} \circ \theta_i$. Twisting $s_{0,V}$ by a 
class in $\rH^1(k,\ov A_V)$ that diagonally maps to $(\alpha_V(s_i))$
we obtain a section $s_V: \Ga_k \to \pi_1^V(X)$ such that $s_V \circ \theta_i$ equals $p_V \circ s_i$ up to conjugation in $\ov A_V$.

Assumption ($\ast$) now implies that the set of such sections $s_V$  is finite. Again by \cite{bourbaki:top1-4} Chapter I \S9.6 Proposition 8, there is a compatible family of sections $(s_V)$ in the projective limit over all possible $\ov V$, which defines a section 
\[
s : \Ga_k \to \varprojlim_{\ov V} \pi_1^V(X) = \pi_1^U(X)
\]
such that $s \circ \theta_i=s_i$ up to conjugation in $\ov A$ by the projective limit argument as in the proof of Proposition~\ref{converse}. This completes the proof of (c).
\enddem

\rem It is worth noting that the additional assumption that 
(\ref{modfund}) has a section allows us to find a genuine section $s$ that "interpolates" the $s_i$ up
to conjugation even in $\ov A$, and not merely a homomorphism or interpolation up to conjugation in $\pi_1 ^U(X)$. 

\medskip

We can prove more under the additional assumption of the collection of fields $(k_i)$ being  arithmetically sufficiently rich in a sense to be made precise as follows. Consider the following property.
\begin{enumerate}
\item[($\ast \ast$)] The union of the conjugates of all the images $\Ga_i \to \Ga_k$ is dense in $\Ga_k$.
\end{enumerate}

\begin{lem} \label{starsfinite}
Property ($\ast \ast$) is inherited by finite extensions $k'/k$ with respect to the set of all composita $k_i\cdot k'$.
\end{lem}
\dem
For any $\sigma \in \Ga_k$ let 
$k'_{i,\sigma}$ be the field extension of $k_i$ associated to the preimage $\Ga'_{i,\sigma} = \theta_i^{-1}(\sigma^{-1}\Ga_{k'}\sigma)$ in $\Ga_i$, namely the compositum $k_i\sigma^{-1}(k')$ in $\kbar_i$ using the fixed embedding $\sigma^{-1}(k') \subset \kbar \subset \kbar_i$ that yields $\theta_i$. The inclusion $k' \subset k'_{i,\sigma}$ induces the map 
$\theta'_{i,\sigma} = \sigma(-)\sigma^{-1} \circ \theta_i : \Ga'_{i,\sigma} \to \Ga_{k'}$.

The union of the conjugates of the images of all $\theta'_{i,\sigma}$ is dense in $\Ga_{k'}$, saying that property ($\ast \ast$) is inherited for finite field extensions $k'/k$ for the new family of fields $(k'_{i,\sigma})$. Indeed, we have to show that 
\[
\bigcup_{i,\sigma} \theta'_{i,\sigma}(\Ga'_{i,\sigma}) = \big(\bigcup_{i,\sigma} \sigma\Ga_i\sigma^{-1}\big) \cap \Ga_{k'}
\]
surjects onto cofinally any finite continuous quotient of $\Ga_{k'}$. It is enough to treat quotients $p_0:\Ga_{k'} \to G_0$ with $\ker(p_0)$ normal in $\Ga_k$, i.e, the map $p_0$ extends to a finite continuous quotient $p:\Ga_k \to G$. Then 
\[
p_0\Big( \big(\bigcup_{i,\sigma} \sigma\Ga_i\sigma^{-1}\Big) \cap \Ga_{k'}\Big) = p\Big(\bigcup_{i,\sigma} \sigma\Ga_i\sigma^{-1}\Big) \cap G_0 = G_0
\]
by property ($\ast \ast$).
\enddem

\begin{lem}  \label{stars}
Property ($\ast \ast$) implies property ($\ast$). 
\end{lem}

\dem
To prove ($\ast$) we consider a finite $k$-group $G$ and $\alpha \in \rH^1(k,G)$. We need to show that the following set is finite:
\[
\Sha_{\alpha} = \{ \beta \in \rH^1(k,G) \ ; \ \theta_i^\ast(\beta) = \theta^\ast_i(\alpha) \in \rH^1(k_i,G) \text{ for all } i \in I \}.
\]
By the technique of twisting, see \cite{cogal} I \S5.4, we may assume that $\alpha$ is the trivial class in $\rH^1(k,G)$. 
Let $k'/k$ be a finite Galois extension that trivialises $G$. With the notation as in Lemma~\ref{starsfinite}, the commutative diagram 
\[
\begin{CD}
\rH^1(k,G) @>{{\rm res}}>> \rH^1(k',G) \cr 
@V{\theta_i^\ast}VV @VV{\theta'^\ast_{i,\sigma}}V \cr 
\rH^1(k_i,G) @>{{\rm res}}>>  \rH^1(k'_{i,\sigma},G) 
\end{CD}
\]
shows that under restriction $\Sha_\alpha$ maps into 
\[
\Sha_{{\rm trivial}} = \{ \chi \in \Hom^\out(\Ga_{k'},G)  \ ; \  \chi \circ \theta'_{i,\sigma} = 1 \text{ for all } i,\sigma \}
\]
which contains only the trivial class due to property ($\ast \ast$) and Lemma~\ref{starsfinite}. Hence, due to the nonabelian inflation--restriction sequence, $\Sha_\alpha$ is contained in $\rH^1(\gal(k'/k),G(k'))$ 
which is a finite set.
\enddem

\begin{prop} \label{denseac}
Under the assumption of ($\ast \ast$)  the properties (a) and (c) are equivalent.
\end{prop}

\dem By Lemma \ref{stars} we also have assumption $(\ast)$.
By Proposition \ref{general} it suffices to show that under assumption (a) the map  $\pi_1^U(X) \to \Ga_k$ admits a section. As (a) trivially implies (a') we may use Proposition \ref{converse} to deduce (b), so that we have found at least a continuous homomorphism $u: \Ga_k \to \pi_1^U(X)$ such that for all $i \in I$ we have $s_i = u \circ \theta_i$ up to conjugation in $\pi_1^U(X)$. Let $\varphi : \Ga_k \to \Ga_k$ 
be the composition of $u$ with the projection $p: \pi_1 ^U(X) \to \Ga_k$. 
To find a section $s_0$ of (\ref{modfund}) and thus to complete the proof of Proposition~\ref{denseac}, it suffices to prove that 
$\varphi$ is bijective because we can then take $s_0=u \circ \varphi^{-1}$. 
We have 
\[
\varphi \circ \theta_i = p \circ (u \circ \theta_i) = p \circ s_i = \theta_i
\]
up to conjugation in $\Ga_k$. Thus for every $\gamma \in  \bigcup_{i} \bigcup_{g \in \Ga_k} g \theta_i(\Ga_i) g^{-1}$ the image $\varphi(\gamma)$ is conjugate to $\gamma$ in $\Ga_k$. By assumption ($\ast \ast$) the set $\bigcup_{i} \bigcup_{g \in \Ga_k} g \theta_i(\Ga_i) g^{-1}$ is dense in $\Ga_k$ so that $\varphi$ preserves every conjugacy class of $\Ga_k$ by continuity and compactness of $\Ga_k$. 
In particular $\varphi$ is injective.

In every finite quotient $\Ga_k \to G$ the image of $\varphi(\Ga_k)$ is a subgroup $H<G$ such that the union of the conjugates of $H$ covers $G$. An old argument that goes back to at least Jordan, 
namely the estimate
\[
|G| = |\bigcup_{g \in G/H} gHg^{-1}| \leq (G:H) \cdot (|H|-1) + 1 = |G| - (G:H) + 1 \leq |G|,
\]
shows that necessarily $H=G$. Thus $\varphi$ is also surjective.
\enddem

\rems (1) The isomorphism $\varphi$ that occurs in the proof of Proposition~\ref{denseac} preserves conjugacy classes of elements, hence is of a very special type which is much studied by group theorists.

(2) In the case of a number field, every  automorphism of $\Ga_k$ is induced by an automorphism of $k$ by a theorem of Neukirch, Uchida and Iwasawa, see \cite{neukirch}, and there are also famous extensions of this result by Pop to function fields. In particular every automorphism of $\Ga_{\bQ}$ is an inner automorphism. 

\begin{prop} \label{densebc}
Under the assumption of ($\ast \ast$)  the properties (b) and (c) are equivalent.
\end{prop}
\dem
Clearly (c) implies (b). For the converse let $u:\Ga_k \to \pi_1^U(X)$ be a homomorphism as in (b), so that there are $\gamma_i \in \pi_1^U(X)$ with $u \circ \theta_i = \gamma_i(-)\gamma_i^{-1} \circ s_i$ for all $i \in I$. With the natural projection $p:\pi_1^U(X) \to \Ga_k$ the proof of Proposition~\ref{denseac} says that the homomorphism $\varphi = p \circ u$ is an isomorphism, so that $s = u \circ \varphi^{-1}$ is a section. With $p(\gamma_i):= \sigma_i$ we compute
\[
\varphi \circ \theta_i = p \circ u \circ \theta_i = p \circ \big(\gamma_i(-)\gamma_i^{-1}\big) \circ s_i =
\sigma_i(-)\sigma_i^{-1} \circ \theta_i,
\]
since $s_i$ is a section and thus $p \circ s_i=\theta_i$. Applying $\varphi^{-1}$ to both sides yields with $\tau_i = \varphi^{-1}(\sigma_i^{-1})$ the equation
\[
\tau_i(-)\tau_i^{-1} \circ \theta_i = \varphi^{-1} \circ \theta_i.
\]
Now the section $s$ interpolates the following
\[
s \circ \theta_i = u \circ \varphi^{-1} \circ \theta_i = u \circ \big(\tau_i(-)\tau_i^{-1}\big) \circ \theta_i
= \big(u(\tau_i)(-)u(\tau_i)^{-1}\big) \circ u \circ \theta_i
\]
\[
=  \big(u(\tau_i)(-)u(\tau_i)^{-1}\big) \circ \big(\gamma_i(-)\gamma_i^{-1}\big) \circ s_i = \big((u(\tau_i)\gamma_i)(-)(u(\tau_i)\gamma_i)^{-1}\big) \circ s_i,
\]
and because of 
\[
p(u(\tau_i)\gamma_i) = \varphi(\tau_i)p(\gamma_i) = \sigma_i^{-1}\sigma_i = 1
\]
we find that $s$ actually satisfies the  stronger interpolation property of (c).
\enddem

\begin{cor} \label{summary}
Under the assumption of ($\ast \ast$)  the properties (a), (a'), (b) and (c) are equivalent to each other and to (a'') together with the existence of  section.
\end{cor}
\dem
This follows immediately by Proposition~\ref{denseac}, Proposition~\ref{densebc}, Lemma~\ref{stars},
Proposition~\ref{converse}, and Proposition~\ref{general}.
\enddem


\section{Results over number fields}

From now on we assume that $k$ is a number field. We still 
consider the exact sequence (\ref{modfund}) 
$$1 \to \ov A \to \pi_1 ^U (X) \to \Ga_k \to 1$$
as above. Let $k_v$ be the completion of $k$ at a place $v$ of $k$. A choice of embeddings $\kbar \to \kbar_v$ of the respective algebraic closures identifies the absolute Galois group $\Ga_v = \gal(\kbar_v/k_v)$ of $k_v$ with the decomposition subgroup of $v$, or more precisely the place of $\kbar$ above $v$ corresponding to the embeding $\kbar \to \kbar_v$.
Hence the map 
$\theta_v : \Ga_v \to \Ga_k$ as defined in the introduction is injective.

\begin{theo} \label{numbers}
Let $S$ be a set of places of $k$ of Dirichlet density $0$, for example a finite set of places. Assume that $X(k_v) \neq \emptyset$ for $v \not \in S$.
For each $v \not \in S$, let $s_v : \Ga_v \to \pi_1 ^U(X)$ be the section map
associated to a $k_v$-rational point $P_v \in X(k_v)$. 

Then the following assertions 
are equivalent.
\begin{enumerate}
\item[(i)] For every finite $k$-group scheme $G$ and every $\alpha \in
\rH^1(X,G)$ such that $\bar \alpha$ has trivial restriction to $\ov U$,
the family $(\alpha(P_v))$ 
belongs to the diagonal image of $\rH^1(k,G)$ in 
$\prod_{v \not \in S} \rH^1(k_v,G)$. 
\item[(i')] For every finite constant $k$-group scheme $G$ and every $\alpha \in
\rH^1(X,G)$ such that $\bar \alpha$ has trivial restriction to $\ov U$,
the family $(\alpha(P_v))$ belongs to the diagonal image of $\rH^1(k,G)$ in 
$\prod_{v \not \in S} \rH^1(k_v,G)$. 
\item[(i'')] There is a section $s_0 : \Ga_k \to \pi_1^U(X)$ and for every finite $k$-group scheme $G$ and every $\alpha \in
\rH^1(X,G)$ such that $\bar \alpha$ restricts trivially to $\ov U$ and the associated $G$-torsor on $X$ is geometrically connected, the family $(\alpha(P_v))$ belongs to the diagonal image of 
$\rH^1(k,G)$ in $\prod_{v \not \in S} \rH^1(k_v,G)$. 
\item[(ii)] There exists a homomorphism $s : \Ga_k \to \pi_1 ^U(X)$ of (\ref{modfund})
such that for each $v \not \in S$, we have $s_v=s \circ \theta_v$ up to
conjugation in $\pi_1^U(X)$.
\item[(iii)] There exists a section $s : \Ga_k \to \pi_1 ^U(X)$ of (\ref{modfund})
such that for each $v \not \in S$, we have $s_v=s \circ \theta_v$ up to
conjugation in $\ov A$, i.e, the sections $s_v$ come from a 
global section $s$.
\end{enumerate}
\end{theo}

\dem  
This is merely a translation of Corollary~\ref{summary} into the number field setting, once we notice that assertion ($\ast \ast$) follows immediately from Chebotarev's density theorem.
\enddem

\rems (1) For $\ov U$ trivial, we have $\pi_1 ^U(X)=\pi_1(X)$ and 
assertion (i)  
means in the language of \cite{stollant}, Definition~5.2,
that the family $(P_v)$ survives every $X$-torsor under 
a finite group scheme $G/k$, while assertion (i') says
 that $(P_v)$ survives every $X$-torsor under
a finite {\it constant} group scheme.
By (i'') this is equivalent to the existence of a section\footnote{Thanks to 
M. Stoll for pointing out the importance of this condition.
} 
together with $(P_v)$ surviving every {\it geometrically connected} $X$-torsor under a finite group scheme.

\smallskip

(2) Even when $X(k) \neq \emptyset$,
it is not sufficient to demand in (i'') that $(P_v)$ survives
every geometrically connected torsor
under a finite and constant group scheme to deduce that $(P_v)$ satisfies
the equivalent properties of Theorem~\ref{numbers}.
Take for example $k=\Q$ and $X$ such that $\pi_1(\ov X)=\mu_3$ with the
corresponding Galois action. Such examples arise among varieties of general
type. Then the only $G$-torsor over $X$ with
$G$ finite constant and $Y$ geometrically connected is $X$ with trivial group $G$.
Nevertheless, there is a geometrically connected torsor $Y \to X$ under
$\mu_3$, and certain families
$(P_v)$ do not survive $Y$, see \cite{dhens}, Remark after Corollary~2.4.

\smallskip

(3) An interesting case is 
when $\ov U$ is the closure of the derived subgroup of 
$\pi_1(\ov X)$,  so that $\ov A=\pi_1(\ov X)/\ov U$ is just 
the abelianized profinite group $\pi_1^{\rm ab}(\ov X)$.
Then the section $s$ in assertion (iii) corresponds to a section 
of the geometrically abelianized fundamental exact sequence
\[
1 \to \pi_1^{\rm ab}(\ov X) \to \pi_1^U(X) \to \Ga_k \to 1.
\]
Then assertion (i) means that the family $(P_v)$ survives every $X$-torsor $Y$ under a finite group scheme $G$ such that $Y \to X$ has \textbf{abelian geometric monodromy}, 
that is: such that the image of the homomorphism $\pi_1(\ov X) \to \ov G$ 
associated to $\ov Y \to \ov X$ is an abelian group.
Similar statements hold for abelian replaced by solvable or nilpotent 
taking for $\ov A$ the maximal prosolvable or pronilpotent quotient of $\pi_1(\ov X)$.

\smallskip

(4) The analogue of Theorem~\ref{numbers} holds with the same proof if 
we replace the family of completions $(k_v)$ by the corresponding 
henselizations $(k_v^h)$, simply because the assertions only depend on the associated sections and the canonical map $\Ga_{k_v} \to \Ga_{k_v^h}$ is an isomorphism.


\section{"Abelian" applications}

Let $k$ be a number field. Denote by $\Omega_k$ the set of all 
places of $k$.  For a smooth and projective $k$-variety $X$
its {\df Brauer--Manin set} is the subset $X(\A_k)^{\br}$ of 
$\prod_{v \in \Omega_k} X(k_v)$ consisting of those adelic points that 
are orthogonal to the Brauer group for the Brauer--Manin 
pairing, cf. \cite{skobook}, II. Chapter 5.

\smallskip 

The following corollary is a consequence of the implication 
(iii) $\Rightarrow$ (i) in Theorem~\ref{numbers}. 
Similar results had already been observed independently (at least)
by J-L.~Colliot-Th\'el\`ene, 
O.~Wittenberg and the second author. 

\begin{cor} \label{fundbm}
Let $X$ be a smooth, projective, geometrically connected curve over a 
number field $k$. Assume that the abelianized fundamental exact sequence
$$1 \to \pi_1^{\rm ab}(\ov X) \to \Pi \to \Ga_k \to 1$$
has a section $s$, such that for each $v \in \Omega_k$ the corresponding 
section $s_v$ is induced  by a $k_v$-point $P_v$ of $X$.
Then $(P_v) \in X(\A_k)^{\br}$.
\end{cor}

\dem We take for $\ov U$ the closure of the derived subgroup of
$\pi_1(\ov X)$ in Theorem~\ref{numbers}. Then $(P_v)$ satisfies condition 
(iii) of Theorem~\ref{numbers}, 
hence by (i) it survives every $X$-torsor with abelian geometric monodromy under a finite $k$-group $G$. In particular, the adelic point $(P_v)$ survives any $X$-torsor under a finite 
abelian group scheme. It remains to apply \cite{stollant}, Corollary~7.3.
\enddem

\rems
(1) Let $X$ be a smooth projective curve of genus $0$. Then the assumption of Corollary~\ref{fundbm} seems vacuous as $\pi_1(\ov X)=1$ and there is a section with no arithmetic content. But we also assume the existence of an adelic point, whence the curve $X$ has $k$-rational points by the classical Hasse local--global principle for quadratic forms. Moreover, any adelic point on $X\cong \bP_k^1$ satisfies the Brauer--Manin obstruction because $\br(k) =\br(\bP^1_k)$.

(2) Let $X$ be a smooth projective curve of genus $1$ as in Corollary \ref{fundbm} with Jacobian $E$.  Then $X$ corresponds to an element $[X]$ in the Tate--Shafarevich group $\Sha(E/k)$. The existence of an adelic point which survives the Brauer--Manin obstruction then implies by \cite{skobook} Theorem~6.2.3 that $[X]$ belongs to the maximal divisible subgroup of $\Sha(E/k)$, which also follows from \cite{dhsza} Proposition~2.1. When $\Sha(E/k)$ is finite as is conjecturally always the case, then the curve $X$ has a $k$-rational point and is actually an elliptic curve $E$.

\begin{theo} \label{biratabel}
Let $X$ be a smooth, projective and geometrically connected curve over a
number field $k$. Assume that the birational fundamental exact sequence
\begin{equation} \label{birat}
1 \to \Ga_{\kbar(X)} \to \Ga_{k(X)} \to \Ga_k \to 1
\end{equation}
has a section. Then $X(\A_k)^{\br} \neq \emptyset$.
If we assume further that
the Jacobian variety
of $X$ has finitely many rational points and finite Tate--Shafarevich
group, then $X(k) \neq \emptyset$.
\end{theo}

A "non-abelian" version of this theorem will be given in the next 
section (Theorem~\ref{biratunobstructedpoint}).

\dem
A section $s:\Ga_k \to \Ga_{k(X)}$ of (\ref{birat}) 
induces\footnote{This would not be clear if we had replaced $k$ by $k_v$ instead of $k_v^h$. Indeed
the existence of a birational section is not a condition that is stable by extension
of scalars; see \cite{ew:birationalabelian}, Remark 3.12(iii).} 
for every place $v$ of $k$ a section $s_v^h$ for the analogous sequence for $k$ replaced by $k_v^h$. We follow an argument by Koenigsmann \cite{koenigsmann:section} Proposition 2.4 (1).

The image of $s_v^h$ defines a field extension $L_v^h/k_v^h(X)$ as the fixed field in the algebraic closure of $k_v^h(X)$. Because the natural maps between absolute Galois groups
\[
\Ga_{L_v^h} \to \Ga_{k_v^h} \leftarrow \Ga_{k_v}
\]
are isomorphisms, the fields $L_v^h$, $k_v^h$ and $k_v$ are $p$-adically closed fields, see \cite{koenigsmann:p-adic} Theorem 4.1, and thus $L_v^h$ is an elementary extension of $k_v^h$,
see \cite{koenigsmann:section} Fact 2.2. In particular, the tautological $L_v^h$ point of $X$ given by
\[
\Spec(L_v^h) \to \Spec k_v^h(X) \to X
\]
implies the existence of a $k_v^h$-point and thus a $k_v$-point on $X$.

The core of the following well-known limit argument goes back at least to Neukirch, and was introduced to anabelian geometry by Nakamura, while Tamagawa emphasized its significance to the section conjecture. We perform the limit argument by applying the above existence result to every connected branched cover $X' \to X$ (necessarily geometrically connected over $k$) with $s(\Ga_k) \subset \Ga_{k(X')} \subset \Ga_{k(X)}$. Thus the projective system $\varprojlim_{X'} X'(k_v)$ over all such $X'$ is a projective system of nonempty compact spaces, and is therefore itself nonempty by \cite{bourbaki:top1-4} Chapter I \S9.6 Proposition 8.

Let $(P'_v)$ with $P'_v \in X'(k_v)$ be an element in the projective limit with lowest stage $P_v \in X(k_v)$. It follows that the composition of the section $s_{P_v} : \Ga_{k_v} \to \Ga_{k_v(X)}$ with the natural projection $\Ga_{k_v(X)} \to \Ga_{k(X)}$ agrees with the $v$-local component $s \circ \theta_v$ for the original section $s$. We may now apply 
Corollary~\ref{fundbm} to the composition
\[
\Ga_k \xrightarrow{s} \Ga_{k(X)} \to \pi_1(X),
\]
which shows that the adelic point $(P_v)$ of $X$ is orthogonal to 
$\br X$ for the Brauer-Manin pairing.

\smallskip

Under the further assumptions that the Jacobian of $X$ has finite Mordell--Weil group and finite Tate--Shafarevich group we now refer to \cite{stollant} 
Corollary~8.1 (a fact that was observed before by Scharaschkin and 
Skorobogatov) to complete the proof of the theorem.

\enddem

\rem In \cite{ew:birationalabelian} Theorem 2.1 H.~Esnault and O.~Wittenberg discuss a geometrically abelian version of Theorem~\ref{biratabel} with the result that an abelian birational section  yields a divisor of degree $1$ on $X$ under the assumption of the Tate-Shafarevich group of the Jacobian of $X$ being finite.

\bigskip

We next describe an application towards the birational version of the section conjecture of Grothendieck's. Recall that for a geometrically connected $k$-variety $X$ a $k$-rational point $a \in X(k)$ describes a $\pi_1(\ov X)$-conjugacy class of sections $s_a$ of (\ref{fund}). In the birational setting the $k$-rational point leads to the following. Define $\hat{\bZ}(1)$ as the
inverse limit (over $n$) of the $\Ga_k$-modules $\mu_n(\kbar)$.
Due to the characteristic zero assumption the decomposition group $D_a$ of $a \in X(k)$ in $\Ga_{k(X)}$ is an extension
\begin{equation} \label{decomp}
1 \to \hat{\bZ}(1) \to D_a \to \Ga_k \to 1
\end{equation}
that splits for example by the choice of a uniformizer $t$ at $a$ and a compatible choice of $n^{th}$ roots  $t^{1/n}$ of $t$. It follows that up to conjugacy by $\hat{\bZ}(1)$, the inertia group at $a$, we have a {\df packet of sections} of (\ref{decomp}) with a free transitive action by the huge uncountable group
\[
\rH^1(k,\hat{\bZ}(1)) = \varprojlim_n k^\ast/(k^\ast)^n.
\]
It can be proven in at least two different ways that the map $D_a \to \Ga_{k(x)}$ maps the $\hat{\bZ}(1)$-conjugacy classes of sections of (\ref{decomp}) injectively into the set of $\Ga_{\kbar(X)}$-conjugacy classes of sections of (\ref{birat}), see for example \cite{koenigsmann:section} Section 1.4, or \cite{stix:cuspidalex} Section 1.3 and Theorem 14+17, or \cite{eh}.

The {\df birational form of the section conjecture} speculates that for a smooth, projective geometrically connected curve the map from $k$-rational points to packets of sections of (\ref{birat}) is bijective and that there are no other sections of (\ref{birat}), see \cite{koenigsmann:section} Section 1.4+5.

The following theorem is a corollary\footnote{Note however that in the
proof of \cite{stollbis} Theorem 9.18, it is not explained 
why the existence of a birational section over $k$ 
implies the same property over $k_v$.}
of Stoll's results \cite{stollbis}
(Corollary 8.6 and Theorem 9.18).   

\begin{theo} \label{biratSC}
Let $X$ be a smooth, projective and geometrically connected curve over $k$. If we assume that there is a nonconstant map $X \to A$ to an abelian variety $A/k$ with finitely many $k$-rational points and
 finite Tate--Shafarevich group, then every section $s$ of the birational fundamental exact sequence
\begin{equation} \label{birat2}
1 \to \Ga_{\kbar(X)} \to \Ga_{k(X)} \to \Ga_k \to 1
\end{equation}
is the section $s_a$ associated to a $k$-rational point $a \in X(k)$. In other words, the birational section conjecture is true for such curves $X/k$.
\end{theo}


\dem
Let $s$ be a section of (\ref{birat2}), and let $X' \to X$ be a finite branched cover, such that upon suitable choices of base points the image of $s$ is contained in $\Ga_{k(X')} \subset \Ga_{k(X)}$.

Then Theorem~\ref{biratabel} shows that $X'(\A_k)^{\br} \neq \emptyset$.
Exploiting the finite map $X' \to X \to A$ we may use \cite{stollant} 
Theorem 8.6, to deduce $X'(k) \not= \emptyset$. Cofinally all such $X'$ will have genus at least $2$
so that $X'(k)$ then is nonempty and finite by Faltings--Mordell \cite{faltings} Satz 7. It follows that $\varprojlim_{X'} X'(k)$, where $X'$ ranges over the system of all $X'$ as above, is nonempty by \cite{bourbaki:top1-4} Chapter I \S9.6 Proposition 8. Let $a \in X(k)$ be the projection to $X(k)$ of an element of $\varprojlim_{X'} X'(k)$, then the image of $s$ is contained in the decomposition subgroup $D_a \subset \Ga_{k(X)}$ and $s$ belongs to the packet of sections associated to the $k$-rational point $a$.

It remains to refer to the literature for the injectivity of the (birational) section conjecture, which was already known to Grothendieck \cite{letter}, see for example \cite{stix:cuspidalex} Appendix B.
\enddem

\rems (1) By \cite{mazur:eisenstein} every Jacobian $J_0(p)$ of the modular curve $X_0(p)$ for $p=11$ or a prime $p \geq 17$ has a quotient, the Eisenstein quotient $\tilde{J}$ \cite{mazur:eisenstein} II (10.4), with finite Mordell--Weil group $\tilde{J}(\bQ)$ and finite Tate--Shafarevich group $\Sha(\tilde{J}/\bQ)$.
Consequently, every smooth, projective geometrically connected curve $X$ over $\bQ$ with a nonconstant map $X \to J_0(p)$ for a prime $p$ as above will satisfy Theorem~\ref{biratSC} and thus the birational section conjecture will hold for such $X$ with $k=\bQ$. Such curves exist for example as smooth iterated hypersurface sections in $J_0(p)$ by hyperplanes of sufficiently high degree due to Bertini's Theorem.

\smallskip

(2) A recent result of Mazur and Rubin, \cite{mazurrubin} Theorem 1.1, guarantees for any algebraic number field $k$ the existence of infinitely many elliptic curves $E/k$ with $E(k) = 0$. As conjecturally $\Sha(E/k)$ is always finite, these elliptic curves and moreover their branched covers $X \to E$ can be used in Theorem \ref{biratSC} to at least conjecturally produce examples of the birational section conjecture over any algebraic number field.


\section{"Non-abelian" applications} 

\begin{theo} \label{genlinear}
Let $X$ be a smooth, projective, geometrically connected curve over a
number field $k$. Let $(P_v)_{v \in \Omega_k}$ be an adelic point of $X$ 
that survives every $X$-torsor under a finite group scheme. Then $(P_v)$ 
survives every $X$-torsor under a linear group scheme. 
\end{theo}

\dem We apply Theorem~\ref{numbers} in the case $\ov U=0$, 
$\ov A=\pi_1(\ov X)$. The hypothesis means that $(P_v)$ satisfies 
condition (i) of this Theorem, hence there is a section 
$s : \Ga_k \to \pi_1(X)$ as in condition (iii). Let $Y \to X$ be 
a geometrically connected torsor under a finite group scheme. 
Using the section $s$, we can lift $(P_v)$ to an adelic point 
$(Q_v)$ on some twisted torsor $Y^{\sigma}$ such that $s$ take values 
in the subgroup $\pi_1(Y^{\sigma})$ of $\pi_1(X)$. This means that
$(Q_v)$ again satisfies condition (iii) of 
Theorem~\ref{numbers}. In particular Corollary~\ref{fundbm} 
implies that $(Q_v) \in Y^{\sigma}(\A_k)^{\br}$. So we have 
proved that for every geometrically connected torsor $Y$ under a 
finite group scheme
$G$, the adelic point $(P_v)$ can be lifted to an adelic point $(Q_v) 
\in Y^{\sigma}(\A_k)^{\br}$ for some twisted torsor $Y^{\sigma}$. 

\smallskip

This still holds if $Y$ is not assumed to be geometrically connected:
indeed the assumption that there exists an adelic point of $X$ 
surviving every $X$-torsor under a finite group scheme implies
(by a result of Stoll, see also \cite{demant}, beginning of 
the proof of Lemma 3) that there exists a geometrically 
connected torsor $Z \to X$ under a finite $k$-group scheme 
$F$, a cocycle $\sigma \in Z^1(k,G)$, and a morphism $F \to  G^{\sigma}$ 
such that $Y^{\sigma}$ is obtained by pushout of the torsor 
$Z$. We conclude with the functoriality of the Brauer--Manin pairing.

\smallskip

It remains to apply the main result of \cite{demant} to finish 
the proof, namely that the \'etale Brauer--Manin obstruction is a priori stronger than the descent obstruction imposed by linear algebraic groups.
\enddem

\rems (1) The previous result does not hold in higher dimension. 
For example there are geometrically rational surfaces $X$, 
in particular $\pi_1(\ov X)=1$, with $X(k) \neq \emptyset$, 
but such that some adelic points $(P_v)$ do not belong to 
$X({\bf A}_k) ^{\br}$. 
For an example with 
an intersection of two quadrics in ${\bf P}^4$ see \cite{ctscras}, p.\ 3, Example a. 
By \cite{skobook}, Theorem 6.1.2 (a), 
such adelic points do not survive the
{\it universal} torsors, which are those torsors under 
the N\'eron-Severi torus of $X$ whose type in the sense of 
Colliot-Th\'el\`ene and Sansuc's descent theory is an isomorphism, see \cite{skobook}, 
Definition~2.3.3.

\smallskip

(2) Let $X$ be a curve of genus at least $2$ such that the
fundamental exact sequence (\ref{fund}) has a section. If we knew
the $p$-adic analogue of Grothendieck's section conjecture,
Theorem~\ref{numbers} and Theorem~\ref{genlinear} 
would yield the existence of an adelic
point $(P_v)$ that survives every torsor under a linear $k$-group
scheme, which is a priori stronger than $(P_v) \in X(\A_k)^{\br}$.
Recall that by
\cite{stollant}, Corollary~8.1.,
the condition $X(\A_k)^{\br} \neq \emptyset$
already implies $X(k) \neq \emptyset$ if the Jacobian variety
of $X$ has finitely many rational points and finite Tate-Shafarevich
group.

\bigskip

The following result is the "non-abelian" version of 
Theorem~\ref{biratabel}.

\begin{theo} \label{biratunobstructedpoint}
Let $X$ be a smooth, projective and geometrically connected curve over a
number field $k$. Assume that the birational fundamental exact sequence
\begin{equation}
1 \to \Ga_{\kbar(X)} \to \Ga_{k(X)} \to \Ga_k \to 1
\end{equation}
has a section. Then $X$ contains an adelic point 
$(P_v)$ that survives every torsor under a linear $k$-group
scheme. 
\end{theo}


\dem We proceed exactly as in the proof of Theorem~\ref{biratabel}, 
except that at the end we apply Theorem~\ref{numbers} instead of
Corollary~\ref{fundbm}, so that we obtain that the adelic point $(P_v)$ of $X$ 
survives every torsor under a finite $k$-group scheme, hence 
every torsor under a linear $k$-group scheme by Theorem~\ref{genlinear}.

\enddem

\bigskip

Let $\Ga_{\kbar(X)} \to \Ga_{\kbar(X)}^p$ be the maximal pro-$p$ quotient of $\Ga_{\kbar(X)}$ and 
\begin{equation} \label{biratpro-p}
1 \to \Ga_{\kbar(X)}^{p} \to \Ga_{k(X)}^{(p)} \to \Ga_k \to 1
\end{equation}
the pushout of (\ref{birat}) by $\Ga_{\kbar(X)} \to \Ga_{\kbar(X)}^p$. With this exact sequence we can prove the following geometrically pro-$p$ version of Theorem~\ref{biratunobstructedpoint}.

\begin{theo} \label{biratunbstructedpro-p}
Let $p$ be a prime number 
and let $k$ be a number field that contains the $p^{th}$ root of unity.
Let $X$ be a smooth, projective and geometrically connected curve over $k$. Assume that the geometrically pro-$p$ birational fundamental exact sequence 
\[
1 \to \Ga_{\kbar(X)}^{p} \to \Ga_{k(X)}^{(p)} \to \Ga_k \to 1
\]
has a section. Then $X$ contains an adelic point 
$(P_v)$ that survives every torsor under a finite $k$-group
scheme with geometric monodromy a finite $p$-group. 
\end{theo}

\dem
We start as in the proof of Theorem~\ref{biratabel}. The local sections $s_v^h : \Ga_{k_v^h} \to \Ga_{k(X)}^{(p)}$ define liftable sections in the sense of 
\cite{popbirat}, hence by  \cite{popbirat} Theorem B 2) we obtain a 
rational point of $X$ in an elementary extension of $k_v^h$. 
We deduce that $X(k_v^h)$ is nonempty and conclude the proof by the same 
limit argument as in the proof of Theorem~\ref{biratabel}. 
\enddem

\merci We thank U.~Goertz for a stimulating question and M.~\c{C}iperiani, M.~Stoll, T.~Szamuely, and O.~Wittenberg for helpful comments. 
This work started when both authors visited the Isaac Newton Institute for 
Mathematical Sciences in Cambridge, whose excellent working conditions and hospitality are gratefully 
aknowledged.

\end{document}